\documentclass{amsart}
\usepackage{bm,url,amsmath,amssymb,amsthm}
\usepackage{verbatim}
\usepackage{tabls}

\newcommand{\Q}{\mathbb{Q}}
\newcommand{\C}{\mathbb{C}}
\newcommand{\R}{\mathbb{R}}
\newcommand{\Z}{\mathbb{Z}}
\newcommand{\h}{\mathbb{H}}
\newcommand{\w}{\mathcal{H}}
\newcommand{\s}{\Lambda}

\newtheorem{theorem}{Theorem}[section]
\newtheorem{theorem*}{Theorem}
\newtheorem{corollary}[theorem]{Corollary}

\theoremstyle{definition}
\newtheorem{definition}[theorem]{Definition}

\numberwithin{equation}{section}

\title{Improved sphere packing lower bounds \linebreak from Hurwitz lattices}

\author{Stephanie Vance}
\address{Department of Chemistry, Computer Science, \& Mathematics\\
Adams State College\\
208 Edgemont Blvd.\\
Alamosa, CO 81102}
\email{vance.stephaniel@gmail.com}

\date{\today}

\begin{document}

\maketitle

\begin{abstract}
In this paper we prove an asymptotic lower bound for the sphere packing density in dimensions divisible by four. This asymptotic lower bound improves on previous asymptotic bounds by a constant factor and improves not just lower bounds for the sphere packing density, but also for the lattice sphere packing density and, in fact, the Hurwitz lattice sphere packing density.

\end{abstract}

\section{Introduction}
The sphere packing density $\Delta_n$ is defined to be the greatest proportion of $\R^n$ (computed as a limit) that can be covered by congruent solid spheres with disjoint interiors. Currently this value is known only for dimensions $n\leq 3$, and the known densities $\Delta_1 = 1$,  $\Delta_2=\pi/\sqrt{12}$, and $\Delta_3=\pi/\sqrt{18}$ are all achieved by lattice sphere packings (i.e., the sphere centers form a lattice).  (See \cite{AW} for a more detailed history of the sphere packing problem and \cite{CS1}, \cite{CS2}, and \cite{Ha} for an overview of the known optimal sphere packings.)

As is the case for dimensions $1$--$3$, in most dimensions many of the densest known sphere packings are constructed from lattices.  (Given a lattice $\s$, a sphere with radius equal to half the length of the shortest non-zero lattice vectors, i.e.,  $\frac{1}{2}||\s||$, is centered at each lattice point.)  The density of an $n$-dimensional lattice sphere packing with sphere center lattice $\s$ is given by the formula  
$$\Delta(\s) = \frac{||\s||^nV_n}{2^n\det(\s)},$$
where $V_n$ denotes the volume of a unit sphere in $\R^n$ and $\det(\s)$ is equal to the volume of a fundamental lattice region.  (Note that $\det(\s)$ is equal to the absolute value of the determinant of any matrix whose columns form a  basis for $\s$.)  
Unfortunately, even with this explicit density formula the lattice sphere packing density $\Delta_{n,L}$ is unknown for all dimensions $n>8$, with the exception of dimension $24$; 
 see \cite{CS2} and \cite{CK} for more details regarding the known optimal lattices.

Despite the absence of a known general solution to either the sphere packing or lattice sphere packing problems, one can still investigate these problems in higher dimensions with the aid of asymptotic bounds for $\Delta_n$ and $\Delta_{n,L}$.  For example, by looking at saturated sphere packings, i.e., sphere packings in which no more spheres can be added unless the interiors are allowed to overlap, one obtains the lower bound $\Delta_n\geq1/2^n$.  Better asymptotic lower bounds have been obtained over the past century by Minkowski in \cite{Mi1}, by Rogers in \cite{Ro}, by Davenport and Rogers in \cite{DR}, and most recently by Ball in 1992.  In \cite{B} Ball proved $\Delta_{n,L}\geq\zeta(n)(n-1)/2^{n-1}$; here $\zeta(n)$ denotes the Riemann Zeta function which converges to $1$ for large values of $n$.    Note that these lower bounds are all proven using lattice sphere packings and hold for both $\Delta_n$ and $\Delta_{n,L}$.  For a more detailed account of these and other asymptotic lower bounds for $\Delta_n$ and $\Delta_{n,L}$ the reader is referred to \cite{GL} and \cite{Z}. 

The purpose of this paper is to prove the lower bound
\begin{equation}\label{the_lb}
\Delta_{4m, L}\geq\frac{3m \zeta(4m)}{2^{4m-3}e(1-e^{-m})},
\end{equation}
which holds for all $m\geq 2$.
As has been the case for previous asymptotic lower bounds proven for $\Delta_n$ and $\Delta_{n,L}$, this bound is non-constructive and there are known low-dimensional sphere packings with densities significantly larger than the bound.  However, it improves on previous lower bounds for $\Delta_n$ and $\Delta_{n,L}$ and is the best lower bound currently known.  Specifically, it improves by a factor of $3/e$ the asymptotic lower bound previously proven by Ball for dimensions $4m$. (Ball's bound for dimension $n$ not divisible by four continues to be the best known asymptotic lower bound for both $\Delta_n$ and $\Delta_{n,L}$.)  The divisibility condition we impose on the dimension is due to the fact that the bound is proven using $4m$-dimensional Hurwitz lattice sphere packings.  These are lattice sphere packings in $\h^m$ in which the lattices of sphere centers are closed under scalar multiplication by the Hurwitz integer ring $\w = \Z[i, j, \frac{1+i+j+k}{2}]$. (Here $\h$ denotes the quaternion skew field 
$$\left\lbrace a+bi+cj+dk:a,b,c,d\in\R\mathrm{\ and\ }i^2=j^2=k^2 = ijk = -1\right\rbrace,$$
 which as a real vector space is isomorphic to $\R^4$.)    Hurwitz lattice sphere packings are a natural choice to use because in low dimensions divisible by four many of the densest known sphere packings have a Hurwitz structure, that is, their lattices of sphere centers are isometric to a lattice in quaternionic Hermitian space that is a module over the Hurwitz integers.  (Note that a quaternionic vector space $M$ is called Hermitian if it is endowed with a Hermitian product, i.e., a map ${h:M\times M\to\h}$ which is linear with respect to the first variable and conjugate linear with respect to the second; the conjugate of a scalar $a+bi+cj+dk\in\h$ is equal to $a-bi-cj-dk$.  Distance in $M$ comes via the trace of $h(\cdot,\cdot)$,  with the inner product of two vectors $x$ and $y$ given by the formula, ${\langle x,y\rangle = \alpha\mathrm{Tr}_{\mathbb{H}/\mathbb{R}}(h(x,y))}$ where $\alpha>0$ is a specified constant and the value of which does not affect packing densities.) The reader is referred to \cite{V} for a table of the densest known lattice sphere packings in dimensions $4m\leq 24$ that have a Hurwitz structure.  

The next section contains several preliminary results on Hurwitz lattices.  In particular, it contains a theorem regarding the existence of a special quaternionic vector space basis contained in a Hurwitz lattice and a Hurwitz analogue of Hlawka's theorem proven in \cite{Hl}.  Both of these theorems are later used in Section 3 and allow us to adapt the proof techniques used by Rogers in \cite{Ro} and prove \eqref{the_lb} using Hurwitz lattice sphere packings. Afterwards in Section 4 we describe how the proof of \eqref{the_lb} can be modified to  obtain similar lower bounds for Eisenstein and Gaussian lattices in $\C^m$; detailed proofs of these bounds are omitted because they do not improve on Ball's lower bound in dimensions $2m$.  Also, in Section 4 we use the Hurwitz analogue of Hlawka's theorem (proven in Section 2) to prove lower bounds for the densities of Hurwitz lattice packings of special convex bodies (not necessarily spherical) in $\h^m$ which satisfy certain symmetry conditions; in particular, they must be invariant under scalar multiplication by the Hurwitz integer units.    These lower bounds can be regarded as a Hurwitz analogue of the Minkowski-Hlawka theorem for $0$-symmetric convex bodies (see \cite[pp.~202--203]{GL}).  

Note that throughout this paper $\h^m$ is identified with $\R^{4m}$ via the map 
$$(a_1+b_1i+c_1j+d_1k,\dots,a_m+b_mi+c_mj+d_mk)\mapsto(a_1,b_1,c_1,d_1,\dots,a_m,b_m,c_m,d_m),$$ and $B_{4m}$ is used in the to denote the closed unit ball in $\h^m$.  Finally, all integrals represent Riemann integration and while we often work with the quaternionic components of vectors in $\h^m$, in order to perform Riemann integration on real valued functions defined on $\h^m$ we identify $\h^m$ with $\R^{4m}$ via the real vector space isomorphism given above.

\section{Special Properties of Hurwitz lattices}

A Hurwitz lattice is a lattice $\s$ in $\h^m$ that has the extra algebraic structure of a Hurwitz module.  In particular, $\s = \w\s = \lbrace \lambda u: \lambda\in\w \mathrm{\ and\ } u\in\s\rbrace$ where $\w = \Z[i, j, \frac{1+i+j+k}{2}]$ is the ring of Hurwitz integers, a maximal order in the rational positive definite quaternion algebra $\left(\frac{-1, -1}{\Q}\right) = \lbrace a+bi+cj+dk\in\h: a,b,c,d\in\Q \rbrace.$  
As an $\w$-module, every Hurwitz lattice $\s$ in $\h^m$ is necessarily free and is generated by a quaternionic basis for $\h^m$.  (This is due to the fact that $\w$ has class number one, a special property not satisfied by maximal orders in general.)  Furthermore we prove below (see Theorem \ref{thm: sv}) that $\s$ also contains, but is not necessarily generated by, a special quaternionic basis for $\h^m$ with prescribed lengths $\mathrm{min}_1(\s),\dots,\mathrm{min}_m(\s)$ defined as follows.   (See \cite{Ma} and \cite{Re} for additional properties of  properties of Hurwitz lattices, and more generally lattices over maximal orders in positive definite $\Q$-algebras.)

\begin{definition}\label{def: min} The \textit{$i^{th}$ quaternionic minimum $\mathrm{min}_i(\s)$} of a Hurwitz lattice $\s$ in $\h^m$ is the smallest $r>0$ such that the closed ball $r  B_{4m}$ contains $i$ $\h$-linearly independent lattice vectors. 
\end{definition}

Observe that the quaternionic minima defined in Definition \ref{def: min} satisfy $\mathrm{min}_1(\s)\leq \dots\leq \mathrm{min}_m(\s)$ and $\mathrm{min}_1(\s) = ||\s||$, the length of the shortest non-zero lattice vectors in $\s$.  Moreover, they are the quaternionic analogue of the successive minima $\lambda_1,\dots,\lambda_n$ of a lattice in $\R^n$ with respect to the closed unit ball (defined by Minkowski in \cite{Mi2}) and can be generalized to any convex body in $\h^m$ invariant under multiplication by the Hurwitz integer unit group. (See also \cite{GL} regarding the successive minima of an lattice with respect to any convex body.)  

The following theorem, Theorem \ref{thm: sv}, can be regarded as a Hurwitz analogue of a  result proven originally by Minkowski for the successive minima of quadratic forms (and hence for lattices) in \cite{Mi2}.

\begin{theorem}\label{thm: sv} Each Hurwitz lattice $\s$ in $\h^m$ contains a quaternionic vector space basis $v_1,\dots,v_m$ such that $||v_i|| = \mathrm{min}_i(\s)$ for $i = 1,\dots, m$.  Moreover, if \nolinebreak{$\det(\s) = 1$}, then there exists a Hurwitz lattice $\widetilde{\s}$ in $\h^m$ with $\det(\widetilde{\s}) = 1$ and $$||\widetilde{\s}|| = \left(\prod_{i = 1}^m \mathrm{min}_i(\s)\right)^{\frac{1}{m}}.$$  
\end{theorem}

\begin{proof}  
Let $\s$ be a Hurwitz lattice in $\h^m$ and suppose that $\lbrace v_1,\dots,v_{k}\rbrace\subseteq\s$ is a maximal set of $\h$-linearly independent vectors in $\s$ satisfying $||v_i|| = \mathrm{min}_i(\s)$; observe that such a set exists since one can choose $v_1$ to be any minimal vector in $\s$.  If $k<m$ then let $\lbrace w_1,\dots,w_{k+1}\rbrace $ be a set of $\h$-linearly independent vectors, each satisfying  $||w_i|| \leq \mathrm{min}_{k+1}(\s)$.  At least one of these vectors, say $w_{k+1}$, does not lie in the span of $v_1,\dots,v_{k}$ and hence $\lbrace v_1,\dots,v_{k},w_{k+1}\rbrace$ is an $\h$-linearly independent subset of $\s$.  The maximality of $\lbrace v_1,\dots,v_{k}\rbrace$ and the definition of $\mathrm{min}_{k+1}(\s)$ require that $||w_{k+1}||<\mathrm{min}_{k+1}(\s)$. So if $i$ is the smallest index satisfying $||w_{k+1}||<\mathrm{min}_i(\s)$ then the set $\lbrace v_1,\dots,v_{i-1},w_{k+1}\rbrace$ is an $\h$-linearly independent subset of $i$ vectors in $\s$, each vector having length less than $\mathrm{min}_{i}(\s)$. (Note that the condition imposed on $i$ is necessary to conclude this due to the possibility that $\mathrm{min}_{k}(\s) = \mathrm{min}_{k+1}(\s)$.) This is a contradiction.  Therefore $k=m$ and the first claim of the theorem holds.

To prove the second claim, assume that $\det(\s) = 1$ and let $\mathcal{B} = \lbrace b_1,\dots,b_m\rbrace$  be an orthonormal basis (with respect to the Hermitian product $h(x,y) = \overline{x}^Ty$), chosen such that $\mathrm{Span}_\h\lbrace b_1,\dots,b_{{k}}\rbrace = \mathrm{Span}_\h\lbrace v_1,\dots,v_{{k}}\rbrace$ for ${{k}} = 1,\dots,m$. (Note that $\mathcal{B}$ can be obtained by applying Gram-Schmidt orthogonalization to the vectors $v_1,\dots,v_m$ using the Hermitian product $h(\cdot,\cdot)$.)  
Writing each vector in $\h^m$ in coordinates with respect $\mathcal{B}$, define  
$T(x) = \left(\frac{x_1}{\mathrm{min}_1(\s)}, \dots, \frac{x_m}{\mathrm{min}_m(\s)}\right)$ so that $T$ is an $\h$-linear transformation.  Then using $T$ and $\alpha = \prod_{i = 1}^m \mathrm{min}_i(\s)$, define ${\widetilde{\s} = \left\lbrace \alpha^{\frac{1}{m}}T(x) : x\in \s\right\rbrace}$ so that $\widetilde{\s}$ is a Hurwitz lattice in $\h^m$ with 
$$ \det(\widetilde{\s}) = (\alpha^{\frac{1}{m}})^{4m} \det(T)  \det(\s) = \alpha^{4} \alpha^{-4} 
=  1.$$ 
For any non-zero vector $\alpha^{\frac{1}{m}}T(x)\in\widetilde{\s}$ one can choose $k$ such that $x_{k} \ne 0$ and $x_{k+1}=\dots=x_m = 0$.  Due to our choice of basis $\mathcal{B}$, $x\notin\mathrm{span}_\h\lbrace v_1,\dots,v_{k-1}\rbrace$, implying that $||x||\geq\mathrm{min}_k(\s)$.  Therefore, 
$$
\left|\left|\alpha^{\frac{1}{m}}T(x)\right|\right|^2  =  \alpha^{\frac{2}{m}}\sum_{i = 1}^k\left|\frac{x_i}{\mathrm{min}_i(\s)}\right|^2
\geq  \alpha^{\frac{2}{m}}\frac{||x||^2}{\mathrm{min}_k(\s)^2}
\geq  \alpha^{\frac{2}{m}}
$$
Now since $||\alpha^{\frac{1}{m}}T(v_1)||^2 = \alpha^{\frac{2}{m}}$ it follows that $||\widetilde{\s}|| = \alpha^{1/m} = \left(\prod_{i = 1}^m \mathrm{min}_i(\s)\right)^{\frac{1}{m}}.$
\end{proof}

The remainder of this section is devoted to proving the following theorem, Theorem \ref{thm: R2}, a Hurwitz analogue of Hlawka's theorem in \cite{Hl}, and as a corollary, a modified version which involves the factor  $1/\zeta(4m)$.  (Recall $\zeta(\cdot)$ denotes the Riemann Zeta function.)  Note that the corollary, Corollary \ref{cor: mobius}, is referred to in the next two sections.  

\begin{theorem}\label{thm: R2}
Suppose $f:\h^m\to\R$ is a bounded, Riemann integrable function with compact support. If $m\geq 2$, then for every $\varepsilon>0$ there exists a Hurwitz lattice $\s$ in $\h^m$ with determinant one satisfying
\begin{equation}\label{eq: R2}
\sum_{u\in \s\backslash\lbrace0\rbrace}f(u) < \int_{\h^m}f(z)\,dz + \varepsilon.
\end{equation}
\end{theorem}

To prove Theorem \ref{thm: R2} we shall borrow the proof techniques used by Davenport and Rogers in \cite{DR} to prove Hlawka's original theorem.  
In particular, we identify $\h^{m-1}$ with the subspace $\lbrace z\in\h^m: z_m = 0\rbrace$ and for $w\in\h^{m-1}$ and $\alpha>0$ we use $(w,\alpha)$ to denote the vector $w + (0,\dots,0,\alpha)\in\h^m$.  If $\s$ is a Hurwitz lattice in $\h^{m-1}$ then for a fixed choice of $\alpha$ we define $\s_w = \lbrace \s + \lambda(w, \alpha) : \lambda\in\w\rbrace$ so that $\s_w$ is a Hurwitz lattice in $\h^m$ with  $\det(\s_w) = \alpha^4\det(\w)\det(\s) = \frac{\alpha^4}{2}\det(\s).$
Then as done by Davenport and Rogers, for a fixed (Hurwitz) lattice $\s$ and $\alpha>0$ we average the sum on the left hand side of \eqref{eq: R2} over all lattices $\s_w$ and prove the existence of a (Hurwitz) lattice satisfying the inequality.

\begin{proof}
First observe that since $f$ is Riemann integrable we can compute $\int_{\h^m} f(z)\ dz$ as the limit
$$\lim_{\alpha\to0^+}\sum_{\lambda\in\alpha\w}\det(\alpha\w)\int_{\h^{m-1}}f(z,\lambda)\,dz = \lim_{\alpha\to0^+}\sum_{\lambda\in\w\backslash\{0\}}\frac{\alpha^4}{2}\int_{\h^{m-1}}f(z,\alpha\lambda)\,dz.$$
Hence we can choose $\alpha>0$ such that 
$$\sum_{\lambda\in\w\backslash\{0\}}\frac{\alpha^4}{2}\int_{\h^{m-1}}f(z,\alpha\lambda)\,dz<\int_{\h^{m}}f(z)\,dz +\varepsilon.$$Furthermore, since $f$ vanishes off a compact set we can choose this $\alpha$ sufficiently small that $$\s = \left(\alpha^{-1/(m-1)}\det(\w)^{m/(4-4m)}\right){\w^{m-1}}$$ 
 is a Hurwitz lattice in $\h^{m-1}$ with the property $f(u,\lambda) = 0$ for all $u\in\s\backslash\{0\}$ and $\lambda\in\h$.  (Note that $\s$ is scaled so that for every $w\in\h^{m-1}$ the Hurwitz lattice 
$$\s_w = \lbrace \s + \lambda(w, \alpha) : \lambda\in\w\rbrace$$ in $\h^m$ has determinant one.)  

We now show that there exists a vector $w\in\h^{m-1}$ such that $\s_w$ satisfies inequality \eqref{eq: R2}.  To do this, we first suppose that $f$ is continuous.  This assumption, combined with the fact that $f$ is periodic modulo $\s$, allows us to average the left hand side of the equation 
\begin{equation}\label{eq:f}
\sum_{u\in\s_w\backslash\lbrace0\rbrace}  f(u) = \sum_{\lambda\in\w\backslash\lbrace0\rbrace}\sum_{v\in\s}f(v+\lambda w,\lambda\alpha)
\end{equation}
over all possible vectors $w\in\h^{m-1}$.  Observe that this average can be 
computed by summing the averages of the inner sums on the right hand side of \eqref{eq:f} over all vectors $w\in\lambda^{-1}G$, where $\lambda\in\w\backslash\lbrace 0\rbrace$ and $G$ is a fundamental region of $\s$. (For each value of $\lambda$, the inner sum is periodic modulo $\lambda^{-1}\Lambda$ as a function of $w$.) 
The average of each inner sum over all $w\in\lambda^{-1}G$ is equal to 
$$
\frac{1}{\mathrm{vol}(\lambda^{-1}G)}\int_{\lambda^{-1}G}\sum_{v\in\s}f(v + \lambda w, \lambda\alpha)\, dw,
$$ 
which, using the change of variable $z = \lambda w$ with $dz = \frac{\mathrm{vol}(G)}{\mathrm{vol}(\lambda^{-1}G)}\,dw$,  becomes 
$$\frac{1}{\mathrm{vol}(G)}\int_{G}\sum_{v\in\s} 
f(v + z,\lambda\alpha)\, dz $$ and hence is equal to 
$$\frac{1}{\det(\s)} \int_{\h^{m-1}}f(z,\lambda\alpha)\,dz.$$ 
Thus the average of the left hand side of \eqref{eq:f} over all $w\in\h^{m-1}$ is equal to $$\frac{\alpha^4}{2}\sum_{\lambda\in\w\backslash\lbrace0\rbrace} \int_{\h^{m-1}}f(z,\lambda\alpha)\,dz,$$
and so there exists at least one vector $w\in\h^{m-1}$ such that  
$$\sum_{u\in\s_w\backslash\lbrace0\rbrace}f(u)\leq \frac{\alpha^4}{2}\sum_{\lambda\in\w\backslash\lbrace0\rbrace}\int_{\h^{m-1}}f(z,\lambda\alpha)\,dz< \int_{\h^m}f(z)\,dz +\varepsilon.$$

Now if $f$ is not continuous, we can approximate $f$ by a continuous function $g$ satisfying the hypothesis of the theorem.  To do this we use the fact that $f$ is Riemann integrable and so the set of points $D$ at which $f$ fails to be continuous has measure zero.  Moreover, since $f$ vanishes off a compact set and $D$ has measure zero, there exists a compact set $A$ and an open set $B$ such that $\mathrm{vol}(A)<\frac{\varepsilon}{8M}$, $\mathrm{vol}(B)<\frac{\varepsilon}{4M}$, and $D\subset A\subset B\subset\h^m$ (see \cite[\S 11]{Mu} regarding the existence of $A$ and $B$).
Hence by Urysohn's lemma there exists a continuous function $\phi:\h^m\to[0,1]$ satisfying $\phi\vert_A = 1$ and $\phi\vert_{\h^m\backslash B} = 0$.  Using this function $\phi$ and $M = \sup_{z\in\h^m}|f(z)|$ we define $$g = \max\lbrace 2M\phi-M,f\rbrace,$$ so that $g$ is a continuous function with $-M\leq f\leq g\leq M$ and $g\vert_{\h^m\backslash B} = f$. Observe that the latter two conditions satisfied by $g$ imply $$ \int_{\h^{m}}g(z)\, dz - \int_{\h^{m}}f(z)\, dz 
\leq 2M  \mathrm{vol}(B)<\frac{\varepsilon}{2}.$$ 
Now since $g$ is a continuous function satisfying the hypothesis of the theorem, we have shown above that there exists a vector $w\in\h^{m-1}$ such that the Hurwitz lattice $\s_w$ satisfies 
$$\sum_{{u\in\s_w\backslash\lbrace0\rbrace}}g(u) < \int_{\h^{m}}g(z)\, dz + \frac{\varepsilon}{2}.$$  
Hence 
$$\sum_{u\in\s_w\backslash\lbrace0\rbrace}f(u)\leq \sum_{u\in\s_w\backslash\lbrace0\rbrace}g(u)
<  \int_{\h^m} g(z)\,dz + \frac{\varepsilon}{2} <  \int_{\h^m} f(z)\,dz + \varepsilon.$$

\end{proof}

Using a M\"obius inversion argument similar to that used to prove the Minkowski-Hlawka theorem (see \cite[p.~202]{GL}) one can modify inequality \eqref{eq: R2} in Theorem \ref{thm: R2} and obtain Corollary \ref{cor: mobius} below.  This corollary is used to produce the $\zeta(4m)$ factor of the lower bound \eqref{the_lb}.  The details of the proof are provided for completeness; however, since $\zeta(4m)$ converges exponentially quickly to $1$ as $m\to\infty$ one may wish to skip the proof of the corollary and move on to Section 3.  (In the proof of Theorem \ref{thm: R3} in the next section, Theorem \ref{thm: R2} can be used instead of Corollary \ref{cor: mobius} if the $\zeta(4m)$ factor is omitted from inequality \eqref{eq: R3}.)

Note that in the statement of the corollary we use $\s'$ to denote the set of primitive lattice vectors in a Hurwitz lattice $\s$, which is defined to be the set of non-zero lattice vectors which are not a positive integer multiple of another lattice vector.

\begin{corollary} \label{cor: mobius}
If $m\geq 2$ and $f:\h^m\to\R$ is a non-negative, bounded Riemann integrable function with compact support, then for every $\varepsilon>0$ there exists a Hurwitz lattice $\s$ in $\h^m$ with determinant one satisfying
\begin{equation}\label{eq: prim}
\sum_{u\in \s'}f(u) < \frac{1}{\zeta(4m)}\int_{\h^m}f(z)\,dz + \varepsilon,
\end{equation}
where $\s'$ is the set of all primitive lattice vectors in $\s$.
\end{corollary}

\begin{proof} Let  $\mu:\Z^+\to\lbrace-1, 0, 1\rbrace$ denote the $\mathrm{M\ddot{o}bius}$ function, which satisfies 
\begin{enumerate}
 \item $\sum_{k|t}\mu(k) = 
  \begin{cases}
  1 &\ \mathrm{if\ \ }t = 1 \cr
  0 &\ \mathrm{if\ \ }t\ne 1\cr 
 \end{cases}$, and  
 \item for every real number $s>1$, $$\frac{1}{\zeta(s)} = \sum_{k = 1}^\infty\frac{\mu(k)}{k^s},$$ where $\zeta$ is the Riemann zeta function defined by $\zeta(s) = \sum_{k = 1}^\infty k^{-s}$.
\end{enumerate}
Then letting $M>0$ such that $0\leq f(z)\leq M$ and  choosing $\delta>0$ such that $M\delta^{4m}V_{4m}<\varepsilon/2$, define a new function $g:\h^m\to\R$ such that  
$$g(z) =   \begin{cases}\sum_{k = 1}^\infty \mu(k)f(kz)&\mathrm{if\ \ } ||z||\geq\delta \cr
  M & \mathrm{if\ \ } ||z||<\delta.\cr 
 \end{cases} $$ 
Note that $g$ satisfies the hypothesis of Theorem \ref{thm: R2} ($g$ is bounded since $f$ is bounded and vanishes off a compact set) and hence there exists a Hurwitz lattice $\s$ in $\h^m$ with $\det(\s) = 1$ such that
\begin{eqnarray*} 
\sum_{{u\in\s}\atop{u\ne 0}} g(u) & < & \int_{\h^m}g(z)\ dz + \frac{\varepsilon}{2}\\ &  = &\int_{\h^m\backslash\delta  B_n} \sum_{k = 1}^\infty \mu(k)f(kz)\ dz + M\delta^{4m}  V_{4m}  + \frac{\varepsilon}{2}\\
& < & \sum_{k = 1}^\infty \frac{\mu(k)}{k^{4m}}\int_{\h^m\backslash\delta  B_n} f(z)\ dz
 + \varepsilon\\
& \leq & \frac{1}{\zeta(4m)}\int_{\h^m} f(z)\ dz
 + \varepsilon.\\
 \end{eqnarray*}

Now for the set of lattice vectors $A = \lbrace u\in\s: ||u||\geq \delta\rbrace $, let $A'$ denote the vectors in $A$ which cannot be written as a positive integer multiple of another vector contained in $A$.  (Note that $\lbrace u\in\s': ||u||\geq \delta\rbrace\subseteq A'$.) Using the properties of the $\mathrm{M\ddot{o}bius}$  function, the fact that $f$ is non-negative and $g(z)\geq f(z)$ whenever $||z||<\delta$, we obtain 
\begin{eqnarray*}
\sum_{u\in\s\backslash\{0\}} g(u) & = & \sum_{{u\in\s}\atop{0<||u||<\delta}} g(u) + \sum_{u\in A'}\sum_{s = 1}^\infty g(su)\\
& \geq & \sum_{{u\in\s}\atop{0<||u||<\delta}} f(u) + \sum_{u\in A'} \sum_{s = 1}^\infty\sum_{k = 1}^\infty \mu(k)f(ksu)\\
& \geq & \sum_{{u\in\s}\atop{0<||u||<\delta}} f(u) + \sum_{u\in A'}\sum_{t = 1}^\infty\sum_{k|t} \mu(k)f(tu)\\
&  =  & \sum_{{u\in\s}\atop{0<||u||<\delta}} f(u) + \sum_{u\in A'}f(u)\\
& \geq &\sum_{u\in\s'} f(u).
\end{eqnarray*}
Combining these with the previous sequence of inequalities yields inequality \eqref{eq: prim}.
\end{proof}

Note that in the previous proof we defined $g(z)$ as a piecewise function to ensure that it is bounded near the origin.  Issues arise at the origin due to the unknown asymptotic behavior of the Mertens function $M(n) =  \sum_{k = 1}^n \mu(k)$; see \cite{KR} for additional information on both the M\"obius and Mertens functions.

\section{A lower bound for the Hurwitz lattice sphere packing density}

In this section we use Hurwitz lattice sphere packings to prove \eqref{the_lb}, the lower bound for the $4m$-dimensional lattice sphere packing density $\Delta_{4m}$ stated in Section~1.  The proof method we use is very similar to that used by Rogers in \cite{Ro}, in which he proved $\Delta_{n,L}\geq {n\zeta(n)}/({2^{n-1}e(1-e^{-n})})$.  In particular, we first use Theorem \ref{thm: sv} and Corollary \ref{cor: mobius} to prove the existence of a Hurwitz lattice with determinant one and whose product of quaternionic minima satisfies a given lower bound.   We then use this particular result with Theorem \ref{thm: sv} and Mahler's compactness theorem to prove the existence of a Hurwitz lattice satisfying \eqref{the_lb}. 

Note that both the statements and proofs of the following two theorems are similar to Theorems $3$ and $4$ in \cite{Ro}, however, using Theorem \ref{thm: sv} and Corollary \ref{cor: mobius} from of the previous section they are modified to take advantage to the extra symmetries of Hurwitz lattices.

\begin{theorem}\label{thm: R3}   
If $m\geq 2$ and $r>0$ satisfy
\begin{equation}\label{eq: R3}
r^{4m}  V_{4m} < \frac{24m \zeta(4m)}{e(1-e^{-m})},
\end{equation}
then there exists a Hurwitz lattice $\s$ with determinant one satisfying 
$$\prod_{i = 1}^m \mathrm{min}_i(\s) > r^m.$$ 
\end{theorem}

\begin{proof}
For $m\geq 2$ and $r>0$ satisfying $r^{4m}  V_{4m} < \frac{24m \zeta(4m)}{e(1-e^{-m})}$, define $\rho:\h^m\to\R$ to be the radial function\footnote{This function $\rho$ is similar to the $\rho$ defined by Rogers in Theorem 3 in his paper \cite{Ro}.}  
$$\rho(z) = 
\begin{cases}
\frac{1}{4} &\ \mathrm{if\ \ }0\leq {||z||}< re^{(1-m)/(4m)},\cr
\frac{1}{4m}-\log\left(\frac{||z||}{r}\right) &\ \mathrm{if\ \ }re^{(1-m)/(4m)}\leq ||z||\leq re^{1/(4m)},\text{\ and\ }\cr
0 &\ \mathrm{if\ \ }re^{1/(4m)}<||z||\cr
\end{cases}
$$
so that $\rho$ satisfies the hypothesis of Corollary \ref{cor: mobius} and (using polar coordinates) 
\small \begin{eqnarray*}
\int_{\h^m}\rho(z)\ dz
& = & r^{4m}  V_{4m} \frac{e^{1-m}}{4} + \int_{re^{(1-m)/(4m)}}^{re^{1/(4m)}}\left(\frac{1}{4m}-\log\left(\frac{\theta}{r}\right)\right)4m\theta^{4m-1}V_{4m}\ d\theta\\
&= & 
 r^{4m}  V_{4m} \left(\frac{e(1-e^{-m})}{4m}\right)\\
& < &6 \zeta(4m).
\end{eqnarray*}
\normalsize
Observe that we can choose $\varepsilon>0$ such that 
$$\frac{1}{\zeta(4m)}\int_{\h^m}\rho(z)\ dz + \varepsilon < 6.$$
and hence by Corollary \ref{cor: mobius} there exists a Hurwitz lattice $\s$ in $\h^m$ with $\det(\s) = 1$ such that    
$$\sum_{u\in\s'}\rho(u)\leq\frac{1}{\zeta(4m)}\int_{\h^m}\rho(z)\ dz + \varepsilon <6.$$  
By Theorem $\ref{thm: sv}$, $\s$ contains a quaternionic basis $\lbrace v_1,\dots, v_m\rbrace$ such that $||v_{{i}}|| = \mathrm{min}_i(\s)$ for ${{i}} = 1,\dots,m$. Observe that the definition of $\mathrm{min}_i(\s)$ implies that each $v_i$ is necessarily a primitive lattice vector and hence the set  
$$A = \lbrace\lambda v_{{i}}: 1\leq {{i}} \leq m \text{\ and\ } \lambda\in\w^\times\rbrace$$
is contained in $\s'$. Now since $\rho$ is a radial function and the Hurwitz integer unit group $\w^\times$ has size $24$,     
$$24\sum_{{{i}} = 1}^m \rho(v_{{{i}}}) = \sum_{u \in A}\rho(u)\leq \sum_{u\in\s'}\rho(u) < 6.$$
Then since $\rho$ is non-negative the above inequalities imply that for $i = 1,\dots, m$\\ $\rho(v_i)<1/4$ and hence   
$$\rho(v_{{{i}}}) \geq \frac{1}{4m}-\log\left(\frac{||v_{{{i}}}||}{r}\right) = \frac{1}{4m} - \log\left(\frac{\mathrm{min}_i(\s)}{r}\right).$$
Therefore
$$\frac{1}{4}-\sum_{{{i}} = 1}^{m}\log\left(\frac{\mathrm{min}_i(\s)}{r}\right) \leq \sum_{{{i}} = 1}^{m}\rho(v_{{{i}}}) < \frac{1}{4},$$ and so $$\sum_{{{i}} = 1}^{m}\log\left(\frac{\mathrm{min}_i(\s)}{r}\right)> 0.$$
Exponentiating both sides of this final inequality we obtain $$\prod_{i = 1}^m \mathrm{min}_i(\s) > r^m.$$ 
\end{proof}

Theorem \ref{thm: R3} is used in the proof of the next theorem, Theorem \ref{thm: R4}, from which the lower bound 
\begin{equation*}
\Delta_{4m, L}\geq\frac{3m \zeta(4m)}{2^{4m-3}e(1-e^{-m})},
\end{equation*}
immediately follows for all $m\geq 2$. (This inequality is \eqref{the_lb} from Section 1.) We wish to emphasize that due to Theorem \ref{thm: R4}, this lower bound holds not just for the $4m$-dimensional sphere packing and lattice sphere packing densities, but also for the $4m$-dimensional Hurwitz lattice sphere packing density.

\begin{theorem}\label{thm: R4}
For every positive integer $m \geq 2$, there exists a Hurwitz lattice sphere packing $\s$ in $\h^m$ with density 
$$\Delta(\s) \geq \frac{3m \zeta(4m)}{2^{4m-3}e(1-e^{-m})}.$$
\end{theorem}

\begin{proof}
This final lower bound follows readily from Theorem \ref{thm: sv}, Theorem \ref{thm: R3}, and Mahler's compactness theorem.  Note that Mahler's compactness theorem guarantees that every sequence of lattices with determinants bounded above, and minimal vector lengths bounded away from zero, has a convergent subsequence. (See \cite[Ch.~3 \S17]{GL} or \cite{Ma} for the topology of the space of $n$-dimensional lattices and for Mahler's compactness theorem, which is referred to in \cite{GL} as the selection theorem of Mahler.) 

Similar to previous results, we shall first construct a sequence of Hurwitz lattices, all with determinant one, and with densities either exceeding or converging from below to $({3m \zeta(4m)})/({2^{4m-3}e(1-e^{-m})}).$  To do this, let $\lbrace r_t\rbrace_{t = 1}^\infty$ be an increasing sequence of positive real numbers such that the sequence $\lbrace r_t^{4m}  V_{4m}\rbrace_{t = 1}^\infty$ converges from below to  
$$\frac{24m \zeta(4m)}{e(1-e^{-m})}.$$
By Theorem \ref{thm: R3} there exists sequence of Hurwitz lattices $\lbrace {\s}_t\rbrace_{t = 1}^\infty$ in $\h^m$, each with  $\det({\s}_t) = 1$ and  
$$\prod_{i = 1}^m \mathrm{min}_i(\s_t) > r_t^m.$$
Theorem \ref{thm: sv} then implies that there exists another sequence of Hurwitz lattices $\lbrace\widetilde{\s}_t\rbrace_{t = 1}^\infty$, each with $\det(\widetilde{\s}_t) = 1$ and $||\widetilde{\s}_t||> r_t$. Either at least one $\widetilde{\s}_t$ has density $$ \Delta(\widetilde{\s}_t) > \frac{3m \zeta(4m)}{2^{4m-3}e(1-e^{-m})}$$ or 
\item $$\lim_{t\to\infty} \Delta(\widetilde{\s}_t) = \left(\frac{1}{2}\right)^{4m} \frac{24m \zeta(4m)}{e(1-e^{-m})} = \frac{3m \zeta(4m)}{2^{4m-3}e(1-e^{-m})}.$$
If only the latter case holds then we use the fact that the space of $4m$-dimensional Hurwitz lattices with unit determinant is a closed subset of the space of all $4m$-dimensional lattices and so by Mahler's compactness theorem there exists a Hurwitz lattice $\s\in\h^m$ with unit determinant and with density  
$$\Delta(\s)= \frac{3m \zeta(4m)}{2^{4m-3}e(1-e^{-m})}.$$
\end{proof}

\section{Remarks}
In addition to the lower bound \eqref{the_lb} proven in the previous section, Corollary \ref{cor: mobius} can also be used to obtain a lower bound for the optimal density of more general Hurwitz packings in $\h^m$ consisting of copies of an  $\w^\times$-invariant convex body $S$ translated by the vectors of a Hurwitz lattice and such that the interiors of the copies of $S$ are disjoint.  (The term $\w^\times$-invariant convex body is used here to describe a compact convex subset of $\h^m$ with non-empty interior and that is invariant under scalar multiplication by the Hurwitz integer unit group $\w^\times$.  Note that an $\w^\times$-invariant convex body in $\h^m$ is the Hurwitz analogue of a $0$-symmetric convex body in $\R^n$, i.e., a convex body invariant under multiplication by the units $\lbrace-1,1\rbrace$ in $\Z$.)  This type of lower bound is similar to the lower bound given in the Minkowski-Hlawka Theorem for $0$-symmetric convex bodies and is proved using similar techniques; see the proof of the Minkowski-Hlawka theorem in \cite{GL} or \cite{Z}.  The extra assumption that the convex body be $\w^\times$-invariant allows one to improve on the lower bound in the Minkowski-Hlawka theorem by a factor of $12$.  

\begin{theorem}\label{cor: ml}
Let $S$ be an $\w^\times$-invariant convex body in $\h^m$.  If $m\geq 2$ then there exists a Hurwitz lattice packing of translated copies of $S$ such that the density of the packing is at least $$\frac{3\zeta(4m)}{2^{4m-3}}.$$  
\end{theorem}

\begin{proof}  
Without loss of generality assume that $\mathrm{vol}(S) = 24\zeta(4m)$ (otherwise replace $S$ by a dilate having volume $24\zeta(4m)$) and for $\varepsilon>0$, choose $\delta\in(0,1)$ such that $\mathrm{vol}(\delta S) = (24-\varepsilon)\zeta(4m)$.  Let $\rho$ denote  the characteristic function of $\delta S$ and let $\s$ be a Hurwitz lattice in $\h^m$ obtained by applying Corollary \ref{cor: mobius} with this $\varepsilon$ and $\rho$ so that  $\det(\s) = 1$ and $$\sum_{u\in\s'} \rho(u) < \frac{1}{\zeta(4m)}\int_{\h^{m}} \rho(z)\ dz + \varepsilon = 24.$$ Observe that if $\delta S$ contains any non-zero vectors of $\s'$ then it must contain at least $24$ since every Hurwitz lattice has at least $24$ minimal vectors and $\delta S$ is $\w^\times$-invariant (the unit group $\h^\times$ has size $24$).  In particular, $\sum_{u\in\s'} \rho(u)$ is either zero or greater than or equal to $24$.  Thus the above inequality implies that $\sum_{u\in\s'} \rho(u) = 0$ and hence $\s$ intersects $\delta S$ only at the origin.

Now since $\varepsilon>0$ is arbitrary, as a consequence of Mahler's compactness theorem, there exists a Hurwitz lattice with determinant one in $\h^m$ that intersects $S$ only at the origin. Such a Hurwitz lattice can be used to obtain a packing of translates of $\frac{1}{2}S$ and the density of this packing is equal to $\frac{3 \zeta(4m)}{2^{4m-3}}.$
\end{proof}  

Note that all of the proofs in this paper can be easily modified to prove results about $2m$-dimensional Gaussian and Eisenstein lattices in $\C^m$, i.e., lattices which are invariant under scalar multiplication by the Gaussian integers $\mathcal{G} = \Z[i]$ or the Eisenstein integers $\mathcal{E} = \Z[\frac{1+i\sqrt{3}}{2}]$ respectively. However the author has chosen to omit the Gaussian and Eisenstein analogues because the lower bounds obtained for Gaussian and Eisenstein lattice sphere packings in dimension $2m$ are respectively $\frac{m\zeta(2m)}{2^{2m-2}}$ and $\frac{3m \zeta(2m)}{2^{2m-1}}$, with neither exceeding Ball's lower bound in \cite{B}.  The reason why the proof method used in this paper works better for the Hurwitz case is because the Hurwitz integer unit group is sufficiently larger than the size of the Gaussian and Eisenstein integer unit groups.  (Note that for the Gaussian and Eisenstein analogues of Theorem \ref{cor: ml} above, the extra assumption that the convex body in $\C^n$ be $\mathcal{G}^\times$ or $\mathcal{E}^\times$-invariant allows one to improve on the lower bound in the Minkowski-Hlawka theorem by a factor of $2$ and $3$ respectively.)

Finally, the techniques used by Ball to prove his asymptotic lower bound in \cite{B} are entirely different from those used in this paper, as well as the methods used by other authors to obtain weaker lower bounds. It would be interesting if one could adapt Ball's method for Hurwitz, Gaussian, and Eisenstein lattices and improve on the asymptotic lower bound given in this paper using these particular lattices.  
However, we do caution the reader that despite the improvements made to the best known asymptotic lower bounds for $\Delta_{4m}$ and $\Delta_{4m,L}$ in this paper, there is a possibility that using only Hurwitz or Eisenstein lattice sphere packings, or more generally any lattice sphere packing, might be limiting due to the extra structure imposed.  In fact, the densest sphere packings in high dimensions may even be disordered, and Torquato and Stillinger have conjectured in \cite{TS} that such packings might provide an exponential increase in density.  However, it seems difficult to analyze disordered sphere packings and, given our current state of knowledge, imposing additional algebraic structure appears to be the most fruitful approach for improving density bounds.

\section{Acknowledgments}  First and foremost, the author thanks her PhD thesis advisor, Henry Cohn, for his many helpful comments and suggestions on preliminary versions of this paper and for introducing her to the sphere packing problem, and more specifically, to asymptotic lower bounds such as that proven by Keith Ball in \cite{B}.  
The author also wishes to thank Microsoft Research for funding a graduate research assistantship as well as travel during early stages of this work.  Additionally, the author wishes to thank the reviewers of this paper for their comments and suggestions regarding future research directions related to this work.

\end{document}